\documentclass{article}
\usepackage[utf8]{inputenc}

\usepackage{natbib}
\setcitestyle{authoryear,open={(},close={)}}
\usepackage{fullpage}
\usepackage{amsmath}
\usepackage{amssymb}
\usepackage{graphics,graphicx}
\usepackage{hyperref}

\usepackage{marginnote}

\newcommand{\D}{\mathcal{D}}

\newcommand{\T}{\mathcal{T}}
\newcommand{\X}{\mathcal{X}}
\renewcommand{\O}{\mathcal{O}}

\newcommand{\G}{\Gamma}

\title{Tegula- exploring a galaxy of two-dimensional periodic tilings}
\author{R\"udiger Zeller \and Olaf Delgado Friedrichs \and Daniel H. Huson}
\date{July 2020}

\begin{document}

\maketitle

\begin{abstract}
Periodic tilings play a role in the decorative arts, in construction and in crystal structures. Combinatorial tiling theory allows the systematic generation, visualization and exploration of such tilings of the plane, sphere and hyperbolic plane, using advanced algorithms and software.
Here we present a ``galaxy'' of tilings that consists of the set of all  2.4 billion
different types of periodic tilings that have Dress complexity up to 24.
We make these available in a database and provide a new program called Tegula that can be used to search and visualize such tilings.

Availability: All tilings and software and are open source and available here:\\
\url{https://ab.inf.uni-tuebingen.de/software/tegula}.
\end{abstract}

\section{Introduction}

Two dimensional periodic tilings play a role in the decorative arts, prominent examples being the euclidean tilings
that cover the Alhambra palace in Granada, Spain, and 
M.C.~Escher's illustrations of euclidean, hyperbolic and spherical tilings involving reptiles, birds and other shapes  \citep{Schattschneider2010}.
Two dimensional euclidean tilings are used in the construction of floors, walls and roofs. Hyperbolic tilings are used in the analysis of minimal surfaces of three-dimensional crystal structures \citep{Ramsden2009,kolbe2018isotopic}.

\begin{figure}
    \hfil
    {
    \begin{tabular}{ccc}
        \includegraphics[height=0.3\textwidth]{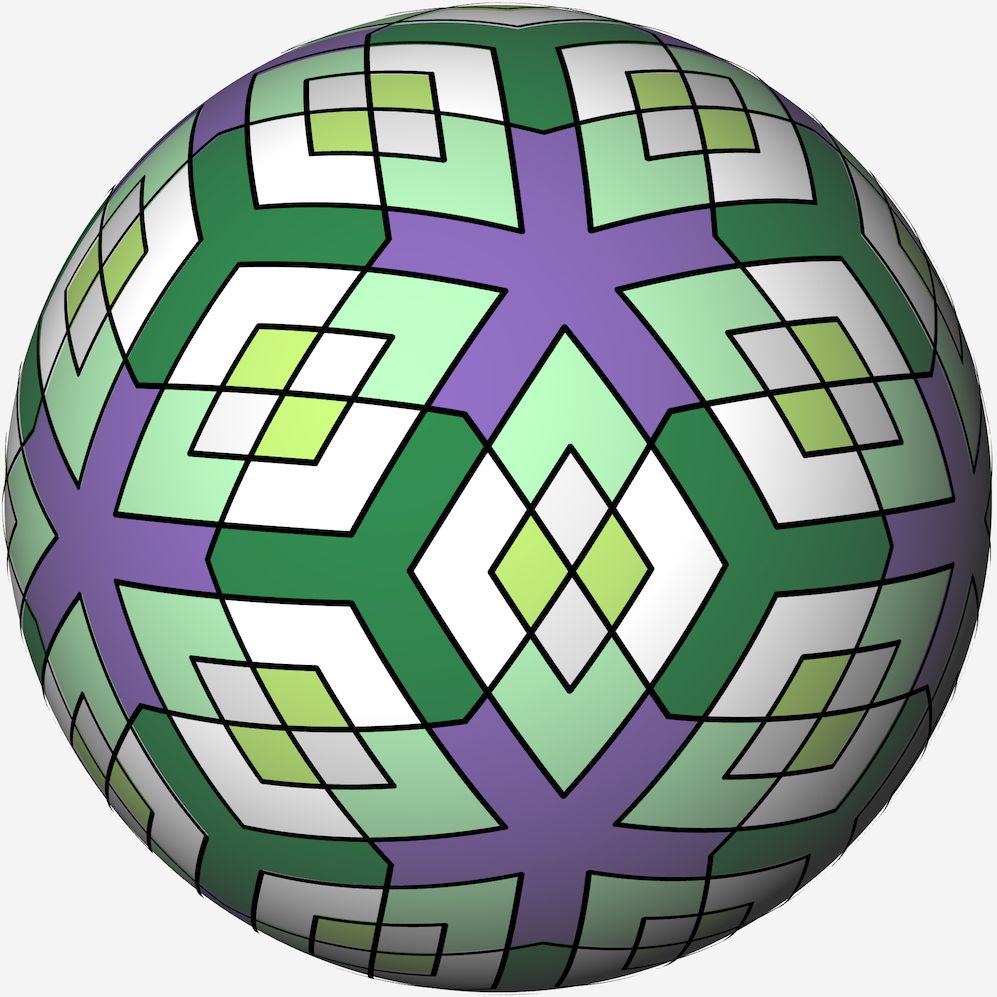} &
        \includegraphics[height=0.3\textwidth]{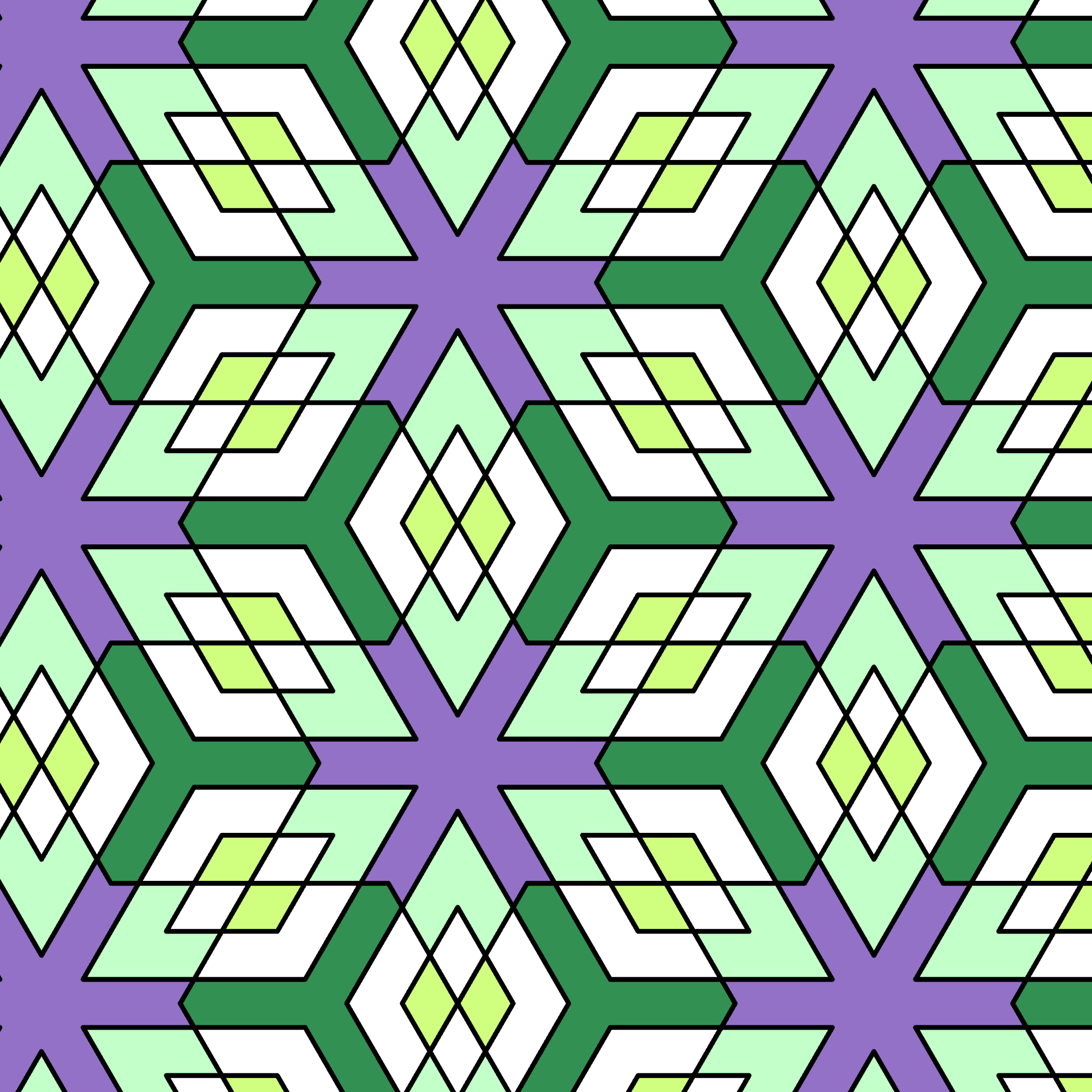} &
    \includegraphics[height=0.3\textwidth]{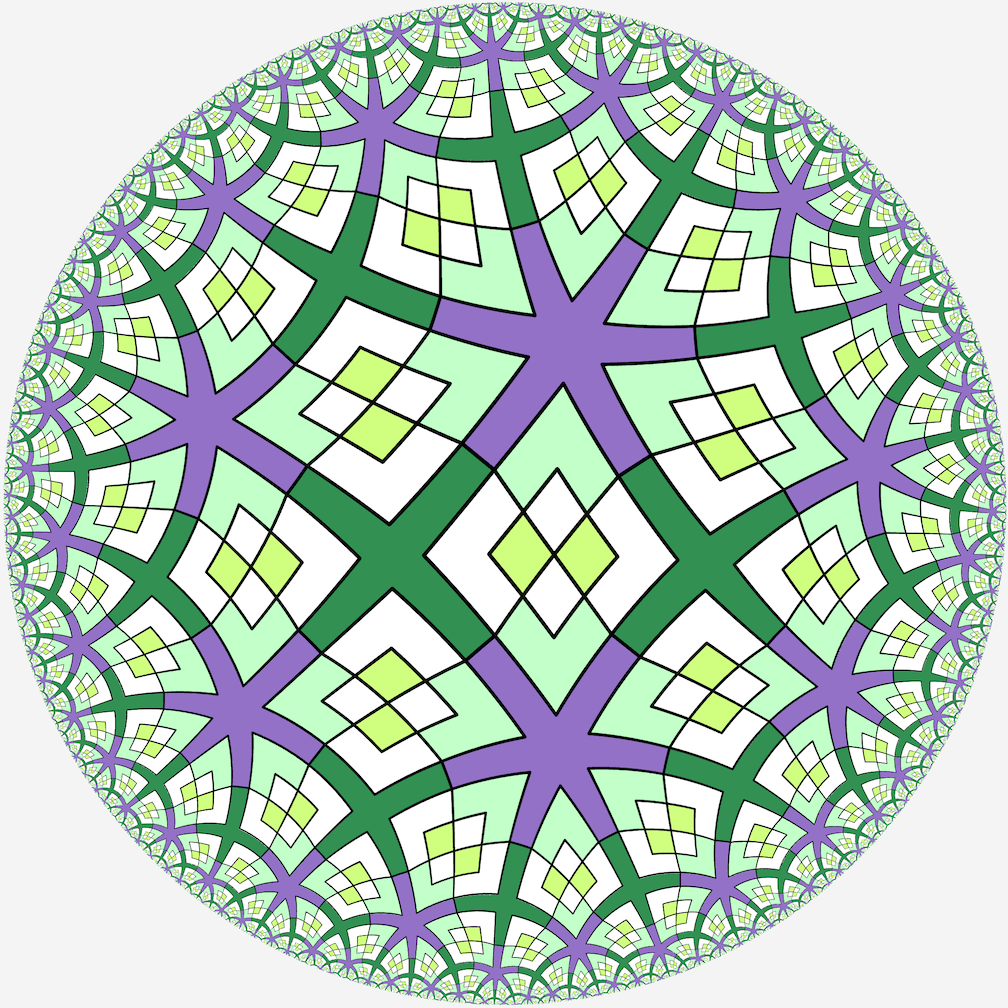} \\
    $\delta=18$, $\G=\mbox{\tt *532}$ & $\delta=18$, $\G=\mbox{\tt *632}$ & $\delta=18$, $\G=\mbox{\tt *642}$\\

    \end{tabular}
    }
    \hfil
    \caption{Periodic tilings of the sphere, plane and hyperbolic plane. Each is labeled by
    its Dress complexity $\delta=\delta(\T,\G)$
    and the orbifold name of the symmetry group $\G$.}
    \label{fig:examples}
\end{figure}

A two-dimensional periodic tiling $(\T,\G)$ of the euclidean plane, sphere or hyperbolic plane, consists of a set of tiles $\T$ 
and a discrete group $\G$ of symmetries of $\T$ with compact fundamental domain, see Figure~\ref{fig:examples}.
Combinatorial tiling theory, based on the encoding of periodic tilings
as ``Delaney-Dress symbols'' \citep{Dress84,Dress87,DressHuson87}, can be used
to systematically enumerate all possible (equivariant) types of two-dimensional
tilings by their curvature and increasing number of equivalence classes of tiles \citep{Huson93a}.

In this paper, we introduce the term  {\em Dress complexity} of a periodic tiling, which is simply the size of the corresponding Delaney-Dress symbol.
We discuss how to systematically enumerate all Delaney-Dress symbols up to
a given Dress complexity, in the case of two-dimensional periodic tilings.
Using this, we have enumerated all two-dimensional periodic tilings of
complexity $\leq 24$. There are $2,395,220,319$ such tilings.
We refer to this collection as a ``galaxy of periodic tilings''
in the title of this paper, because, first, the number of tilings is very big (although not as large as the number of stars in a typical galaxy), and second, when viewing these tilings, the impression is that many look very similar to each other, much like  stars in the sky.

Each such tiling is represented by its Delaney-Dress symbol and we provide
these in a SQLITE database.
We provide a new program called {\em Tegula} that allows the user to explore and
query the database, and to visualize the corresponding tilings in all three geometries.
Tegula and the database of periodic tilings are open source and freely available.

%%%%%%%%%%%%%%%%%%%%%%%%%%%%%%%%%%%%%%%%%
\section{Conway's orbifold notation}

In this section, we briefly recall results on the classification of surfaces \citep{SeifertThrefall1934,ZIPProof} and their Conway names \citep{ConwayHuson99}. Any orientable, closed, connected surface is homeomorphic to either the sphere, denoted by {\tt 1},
or a sphere with $h>0$ handles attached, denoted by
$$\underbrace{\mbox{\tt o o \dots~o}}_h.$$
Any non-orientable, closed, connected surface is homeomorphic to a sphere with $k\geq 1$ crosscaps attached, denoted
by
$$\underbrace{\mbox{\tt x x \dots~x}}_k.$$
The surface {\tt x} is a project plane, the surface {\tt xx} is a Klein bottle
and the surface {\tt xxx} is called Dyck's surface.

Note that the classification of closed surfaces does not mention combining both handles and crosscaps. This is because, if a crosscap is present, then any given handle can be replaced by two crosscaps \citep{Dyck1888}.

A connected surface with boundary is obtained from a closed connected surface by removing $k$ disks from the interior
of the surface.  In Conway's notation, a sphere with $b>0$ boundary components is written as
$$\underbrace{\mbox{\tt * * \dots *}}_b,$$
a sphere with $h>0$ handles and $b>0$ boundary components is written as
$$\underbrace{\mbox{\tt o o \dots~o}}_h\underbrace{\mbox{\tt * * \dots~*}}_b,$$
and
a sphere with $k>0$ crosscaps and $b>0$ boundary components is written as 
$$\underbrace{\mbox{\tt * * \dots~*}}_b\underbrace{\mbox{\tt x x \dots~x}}_k.$$
The surface {\tt *x} is a M\"obius strip.

\begin{figure}
    \centering
    \begin{tabular}{ccc}
    \includegraphics[width=0.3\textwidth]{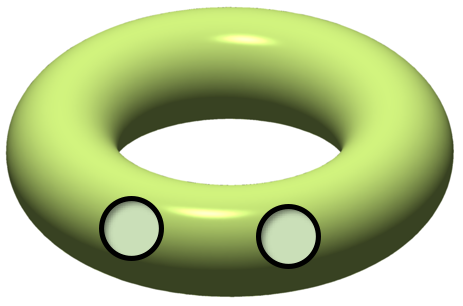} & &
    \includegraphics[width=0.3\textwidth]{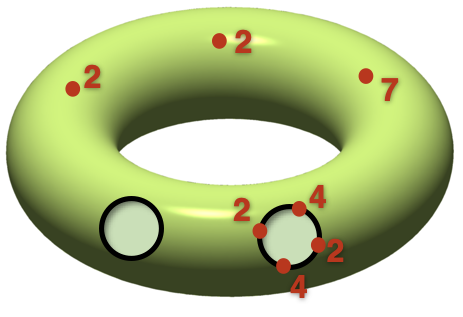}\\
    Surface {\tt o**} & & Orbifold {\tt o227**2424}\\
    \end{tabular}
    \caption{On the left we show the orientable surface {\tt h**}.
    On the right we show the orbifold obtained by adding three cones and four corners.
    (Torus image: Oleg Alexandrov, Wikipedia.)
    }
    \label{fig:orbifold}
\end{figure}

For the purposes of this paper, a {\em two-dimensional orbifold} \citep{Thurston80,ConwayHuson99} consists of a connected surface $S$,
either orientable or non-orientable, with boundary or without, together with a finite set of
points $P=\{p_1,\dots,p_t\}$ in $S$, where each such point $p_i$ is labeled with an integer $v_i\geq 2$ that we call its {\em order}.
Any such point is called a {\em cone}, if it is contained in the interior of $S$, or
a {\em corner},  if it is contained in the boundary of $S$.
For example, in Figure~\ref{fig:orbifold} we depict an orbifold obtained by adding three cones and four corners to the surface {\tt o**}. Note that the set of added points can be empty, so any surface, such
as {\tt o**}, is also an orbifold.

Conway's notation, which we introduced above for the naming of surfaces, also covers orbifolds and has
this form
 \begin{equation}\label{eqn:orbifold}
\underbrace{\mbox{\tt o o \dots~o}}_{h\mbox{~handles}}~
\underbrace{\mbox{\it A B C \dots}}_{\mbox{cones}}~
\underbrace{
*\underbrace{\mbox{\it a b c \dots}}_{\mbox{corners}}~
*\underbrace{\mbox{\it r p q \dots}}_{\mbox{corners}}~
* \dots~}_{b~\mbox{boundary components}}
\underbrace{\mbox{\tt x x \dots~x}}_{k~\mbox{crosscaps}}.
 \end{equation}

Note that the cone degrees $A, B, C, \dots$ are unordered, whereas the corner degrees associated with
any given boundary component have a cyclic ordering given by the order in which they are encountered along
boundary component. One can flip the direction in which corners of an individual boundary component are listed
if there are no other boundary components with corners, or if the surface is non-oriented.

%%%%%%%%%%%%%%%%%%%%%%%%%%%%%%%%%%%%%%%%%
\section{Equivariant tilings}

Through this paper, let $\X$ be one of the three two-dimensional geometries, namely either the
sphere $\mathbb{S}^2$, the euclidean plane $\mathbb{E}^2$, or the hyperbolic plane  $\mathbb{H}^2$.

We  use $(\T,\G)$ to denote a {\em equivariant tiling} of $\X$, defined in the usual way \citep{DressHuson87}. That is, any such tiling $(\T,\G)$ consists
of a  set of tiles  $\T$ that is ``well-behaved'' (i.e.\ a locally finite tiling whose tiles have simply-connected interiors), and a group $\G$ of isometries of $\X$ that map
the set of tiles $\T$ onto itself. 
We emphasize that the word {\em equivariant} indicates that
the symmetry group is prescribed and thus may be only a subgroup
of the group of all automorphisms of the tiling.

Such a tiling $(\T,\G)$ is called {\em periodic} tiling,
if its symmetry group $\G$ is a discrete group with a compact fundamental domain.
Examples of such  tilings are shown in Figure~\ref{fig:examples}.
%\marginnote{Olaf: do we have to say `discrete'? --- Not if we insist that tiles have vertices, strictly speaking, but I would say it anyway.}

\begin{figure}
    \centering
    \begin{tabular}{ccc}
    \begin{tabular}[c]{c}
    \includegraphics[width=0.3\textwidth]{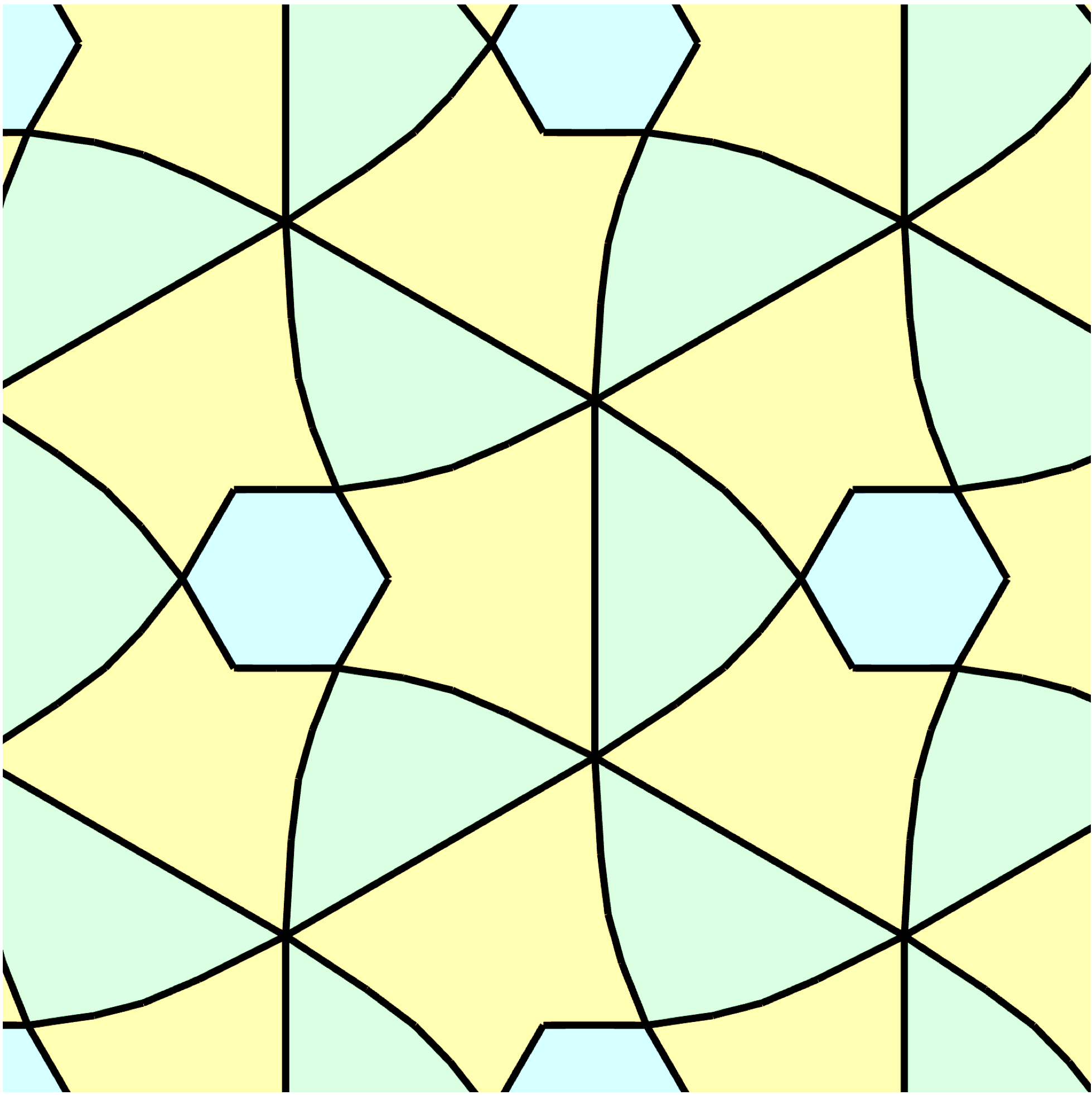}
        \end{tabular}
 &
    \begin{tabular}[c]{c}
    \includegraphics[width=0.3\textwidth]{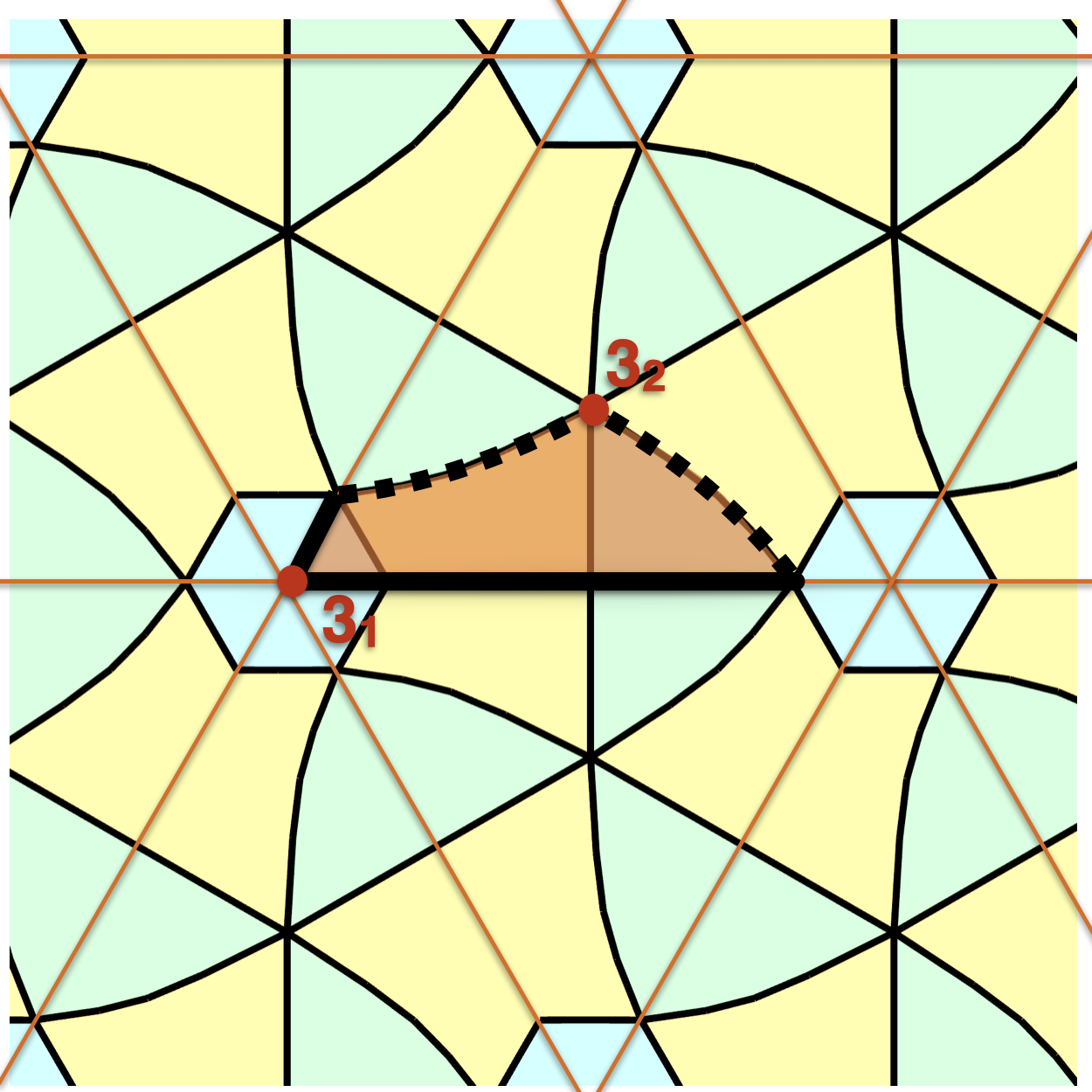}
        \end{tabular}
 &
    \begin{tabular}[c]{c}
     \includegraphics[width=0.15\textwidth]{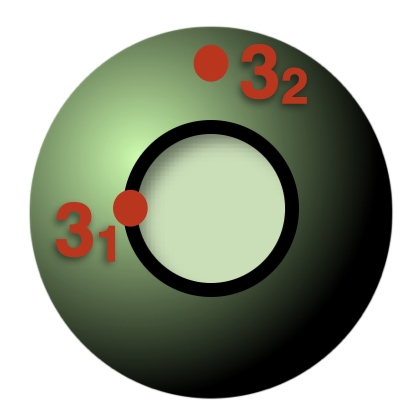}\\
     \\
     {\tt 3*3}\\
     \end{tabular}
     \\
       (a) Periodic tiling $(\T,\G)$ & (b) Fundamental domain
       \& symmetries & (c)  Orbifold\\
    \end{tabular}
    \caption{(a) A periodic tilings of the plane.
    (b) Here we highlight a fundamental domain. Reflectional axes are shown as thin lines.
    The boundary of the fundamental domain that gives rise to the boundary of the
    orbifold is shown as a solid thick line, whereas the two dotted thick lines are identified
    with each other. There are two rotational centers on the boundary of the fundamental domain, labeled $3_1$ and $3_2$, which give rise to a corner and cone, respectively.
    (c) The corresponding orbifold and orbifold name.
        (Sphere image: Darkdadaah, Wikipedia.)
}
    \label{fig:group}
\end{figure}

Let $\G$ be the symmetry group of a two-dimensional
periodic tiling.
The {\em orbifold} $\O(\G)$ of such a group is ``the surface divided by the group'', that is, the orbit-manifold given by
the quotient topological space whose points are the orbits under the group \citep{ConwayHuson99}. Here, reflections are mapped onto boundary segments,
rotational centers are mapped onto cones,
dihedral centers are mapped onto corners, and glide-reflections give rise to crosscaps. The order of such a point is given by
the largest order of any rotation in the symmetry group that fixes the center.
This is illustrated in Figure~\ref{fig:group}, and in Figure~\ref{fig:examples} we
list the orbifold names for the displayed tilings.
Given the drawing of a periodic tiling, or other periodic pattern,
one can easily
determine the orbifold name for corresponding symmetry group, as
discussed in \citep{ConwayHuson99}.

Do all orbifold names correspond to symmetry groups?
Any two-dimensional orbifold with name $\O$, of
the form shown above in Equation~\ref{eqn:orbifold}, can be obtained
as either
\begin{enumerate}
    \item $\mathbb{S}^2~/$ an orthogonal group,
        \item $\mathbb{E}^2~/$ a crystallographic group, or
        \item $\mathbb{H}^2~/$ a ``non-euclidean crystallographic group'',
\end{enumerate}
except for the ``bad orbifolds'' with names
{\tt p}, {\tt pq}, {\tt *p} and  {\tt *pq} with
$p,q \geq 2$ and $p\neq q$, see \citep{ConwayHuson99}.

%%%%%%%%%%%%%%%%%%%%%%%%%%%%%%%%%%%%%%%%%%
\section{Combinatorial tiling theory}

In combinatorial tiling theory, every periodic tiling 
$(\T,\G)$ is represented by a Delaney-Dress symbol $(\D,m)$,  
defined as a finite set $\D$, together with
the action of a certain free group $\Sigma$, together with 
with maps $m_{01}, m_{12}, m_{02}:D \to \mathbb{N}$, fulfilling
certain conditions, see \citep{Dress84,DressHuson87}.

A key result is that the Delaney-Dress symbol describes a periodic tiling up to equivariant equivalence. In more detail, two 
periodic tilings $(\T,\G)$ and $(\T',\G')$ are
equivariantly equivalent, if and only if their corresponding
Delaney-Dress symbols $(\D,m)$ and $(\D',m')$ are isomorphic \citep{Dress84,Dress87}.

Based on this, all two-dimensional periodic tilings can
be systematically enumerated \citep{Huson93a}.
Delaney-Dress symbols can be assigned to higher-dimensional tilings,
and have been used to address classification problems for
three-dimensional euclidean tilings \citep{Molnar1997,DelgadoHuson99a,delgado2005isohedral,dutour2010space} and
as a useful data-structure in the context of developing the
system of orbifold names for three-dimensional euclidean space groups \citep{DelgadoHuson96,ConwayDelgadoHusonThurston2001}.

\begin{figure}
    \centering
    \begin{tabular}{cc}
     \begin{tabular}[c]{c}
    \includegraphics[width=0.4\textwidth]{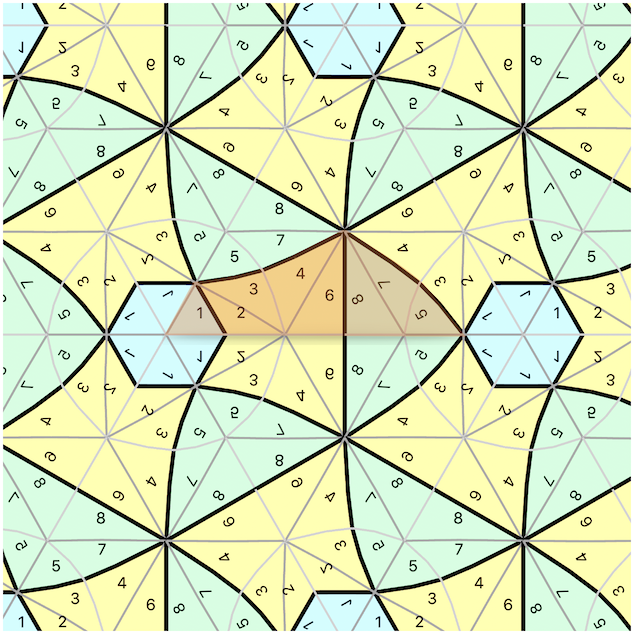} \\
        \end{tabular}
        &
 \begin{tabular}[c]{c}
     \includegraphics[width=0.4\textwidth]{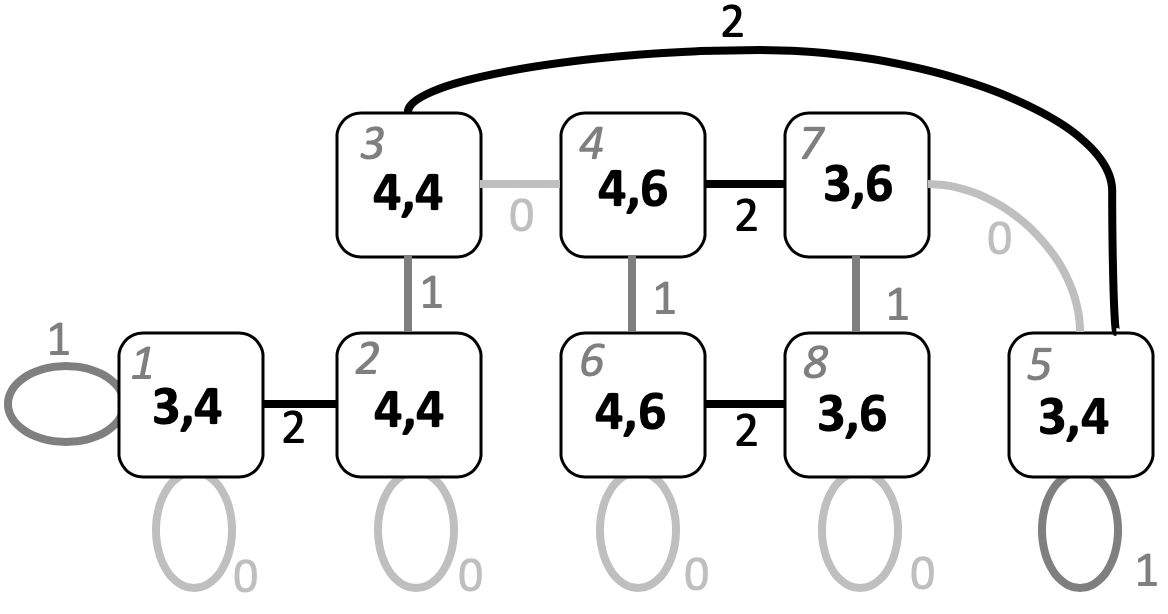}\\
      \end{tabular}
      \\
      (a) Periodic tiling with chamber system & (b) Delaney-Dress symbol\\
  \end{tabular}
    \caption{(a) A periodic tiling, triangulated into chambers,
    all symmetry-equivalent chambers labeled with the same number 1--8.
    The highlighted fundamental domain contains exactly one chamber for each of the eight numbers. (b) The associated
    Delaney-Dress symbol has eight corresponding nodes (labeled
    here 1--8).
    There are three types of edges, labeled 0--2,
    indicating neighbor relationships between chambers.
    Each node is labeled with two bold numbers, indicating
    the degree of the containing tile and the degree of the
    incident vertex, respectively.
    }
    \label{fig:delaney}
\end{figure}

Rather than repeat the details of the definition of a Delaney-Dress symbol here,
we illustrate its construction using an example. 
Consider the periodic tiling shown in Figure~\ref{fig:group}(a).
To construct its Delaney-Dress symbol, start by triangulating the
tiling using a barycentric subdivision, as shown in Figure~\ref{fig:delaney}(a).
Note that each triangle corresponds
to a flag $(v,e,t)$ consisting of a vertex $v$ contained in
an edge $e$, which is contained in a tile $t$. Every triangle
has exactly three neighbors, which we call its $0$-, $1$-
$2$-neighbor, whose flags differs only in their $0$-, $1$- or $2$-component, respectively.

The second step is to partition the set of triangles into
equivalence classes, considering any two triangles to be equivalent,
if their exists an symmetry of the tiling that maps the one triangle onto the other. In this example we obtain eight such classes and label them 1--8.

These eight equivalence classes define the Delaney-Dress set $\D$,
which are represented by nodes in the graph shown
in Figure~\ref{fig:delaney}(b). 
Two such nodes are connected by an edge with label $i$, 
if the $i$-neighbor of any triangle in the one equivalence class
is contained in the other equivalence class.
For example, nodes $1$ and $2$ are connected by an edge labeled $2$,
because triangles of neighboring triangles in classes $1$ and $2$
are incident to the same node and edge, but are contained in different tiles.

The final step is to label each node $D$ by two numbers, $p,q$;
these record the tile-degree (number of edges of the tile) and vertex-degree associated with the given equivalence class of triangles. More formally, the two numbers are denoted by
$m_{01}(D)$ and $m_{12}(D)$.
For example, node $1$ is labeled $3,4$, because all triangles
in equivalence class $1$ are contained in a tile of degree $3$ and are incident to a vertex of degree $4$, whereas node $4$ is labeled $4,6$, because
the corresponding triangles are contained in tiles of degree $4$
and are incident to vertices of degree $6$.

So, we can view a two-dimensional Delaney-Dress symbol as a {\em Delaney-Dress graph}
that is, a connected (multi-) graph in which each node is incident to exactly one edge of each color $0$, $1$, and $2$, together with a labelling of its
nodes by two maps $m_{01}$ and $m_{12}$, fulfilling certain conditions.
 
Let an {\em $i,j$-component} $Z$ be set of nodes in $Z\subseteq \D$ that
is connected by edges of colors $i$ and $j$, with $0\leq i<j\leq 2$.
By the properties of an edge coloring,
$Z$ will always be a cycle or a chain, and
we define the $i,j$-length of $Z$ to be $|Z|$, in the former case,
and $2|Z|$, in the latter. For any node $D\in \D$,
we define $r_{ij}(D)$ to be the $i,j$-length of the $i,j$-component containing $D$.

A natural interpretation of the Delaney symbol is as a triangulation
of the associated orbifold, where each triangle is the image
of one equivalence class of triangles of the original tiling.
By construction, any cone or corner of the orbifold
will lie on a vertex of this triangulation and will be 
surrounded by triangles that belong to the same $i,j$-component
$Z$, for some choice of $0\leq i < j\leq 2$.
We call this a {\em $k$-vertex}, with $k$ such that $\{i,j,k\}=\{0,1,2\}$.
For all nodes $D\in Z$, we define  $v_{ij}(D)$ to be the
order of the associated cone or corner, that is, the highest
order of any rotation in the symmetry group about the vertex.
For each $i,j$-orbit $Z$ whose triangles are not incident to a cone or corner, we set $v_{ij}(D)=1$ for all $D\in Z$.

Note that we have $m_{ij}(D)=v_{ij}(D)r_{ij}(D)$, linking
combinatorial features of the tiling, such as vertex degrees, etc,
with rotational degrees in the symmetry group.
In particular, all equivalence classes of rotational centers
and dihedral centers of the symmetry group can be obtained by
from the Delaney symbol of a tiling by
enumerating all $i,j$-components $Z$ in $\D$ for which
$v_{ij}(D)=\frac{m_{ij}(D)}{r_{ij}(D)}>1$ holds for
$D\in Z$.

The number of symmetry-equivalence classes
of vertices, edges and tiles in a periodic tiling is given by
the number of $1,2$-, $0,2$- and $0,1$-components in its Delaney symbol, respectively.

Other properties of a tiling require more involved analysis
of the corresponding Delaney symbol, such
as the Euler characteristic, curvature, geometry and
the corresponding orbifold name \citep{BalkeHuson94a}.
For example, the curvature is given by the following calculation:
\[
    {\mathcal K}(\D,m)=\sum_{D\in\D}\left(
    \frac{1}{m_{01}(D)}+\frac{1}{m_{12}(D)}-\frac{1}{2}\right),
\]
and this, in turn, defines the geometry associated with the tiling,
namely spherical, euclidean or hyperbolic, depending
on whether the curvature is positive, 0 or negative, respectively.

A more difficult tiling property to obtain from an analysis of the corresponding Delaney symbol is whether the tiling is {\em pseudo convex}, that is, whether the intersection of any two tiles is always either empty or simply connected.

The size of the Delaney-Dress symbol $(\D,m)$ is an important
invariant for the corresponding tiling $(\T,\G)$, albeit
the most simplest property to obtain, and
we propose to call this the {\em Dress complexity} of the tiling,
denoted by $\delta(\T,\G)$.

%%%%%%%%%%%%%%%%%%%%%%%%%%%%%%%%%%%%%%%%%%%%%%%%%%%%%%%%%%%%%%%%%
\section{Enumeration}

A main goal of this paper is to enumerate all periodic tilings
of low Dress complexity.

\begin{figure}
    \centering
    \includegraphics[width=0.3\textwidth]{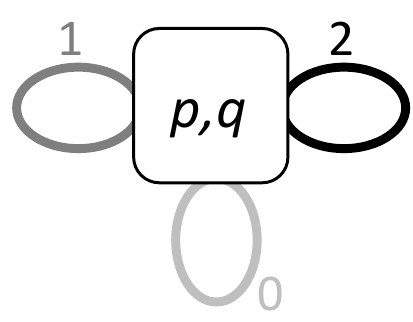}
    \caption{Any Delaney-Dress symbol $(\D,m)$ of Dress complexity $1$
    consists of a single node $D$, three self-edges of colors $0,1,2$,  and two number $p\geq 3$ and $q \geq 3$.}

    \label{fig:symbol1}
\end{figure}

We first start with Dress complexity $1$, that is, Delaney-Dress
symbols of size one, as displayed in Figure~\ref{fig:symbol1}.
In this case, the curvature is given by 
\[{\mathcal K}(\D,m)=\frac{1}{p}+\frac{1}{q}-\frac{1}{2}.\]
For $p=3$ and $q=3,4,5$, this value is positive, and
thus the corresponding tilings are spherical.
The same is true for $p=3,4,5$ and $q=3$.
For $p=3$ and $q=6$, or $p=6$ and $q=3$, or $p=q=4$, the curvature is
$0$ and thus the corresponding tilings are euclidean.
In all other cases, for example $p=4$ and $q=5$, 
the curvature is negative and thus the corresponding tiling is hyperbolic.

If we allow tiles to be digons, that is, to have only two edges, 
then $p=2$ and for any value of $q\geq 3$ the curvature is positive
and so all such tilings of Dress complexity $1$ are tilings of the sphere.
To reduce the number of resulting classes, in this paper
we only enumerate tilings for which all tiles have at least 3 edges,
contrast to some of our previous work \citep{DelgadoHusonZamorzaeva92,Huson93a}.

Already for Dress complexity $1$ we encounter infinite families
of non-equivalent periodic tilings.
To address this, we say that a periodic tiling $(\T,\G)$ 
is {\em geometry minimal}, if one of the three cases hold:
\begin{enumerate}
    \item the tiling is spherical and either the corresponding orbifold is one of $\tt 532$ and $\tt *532$, or all rotational degrees are $\leq 4$,
    \item the tiling is euclidean, or
    \item the tiling is hyperbolic and one can't reduce
     the rotational order of any tile, or vertex,
    without changing the sign of the curvature of the symmetry group,
    or without reducing the degree of the tile, or the vertex, respectively,
    to below $3$.
\end{enumerate}
This property is easily inferred from the corresponding
Delaney-Dress symbol. For a spherical tiling, determine whether
the value of $v_{ij}$ is $\leq 5$ on all $0,1$- and $1,2$-orbits.
For a hyperbolic tiling, reduce the value of $v_{ij}$
on each $0,1$- and $1,2$-orbit in turn, and check whether the 
modified Delaney-Dress symbol has negative curvature and that
the resulting value for $m_{ij}$ is $\geq 3$.

\begin{figure}
    \centering
    \includegraphics[width=\textwidth]{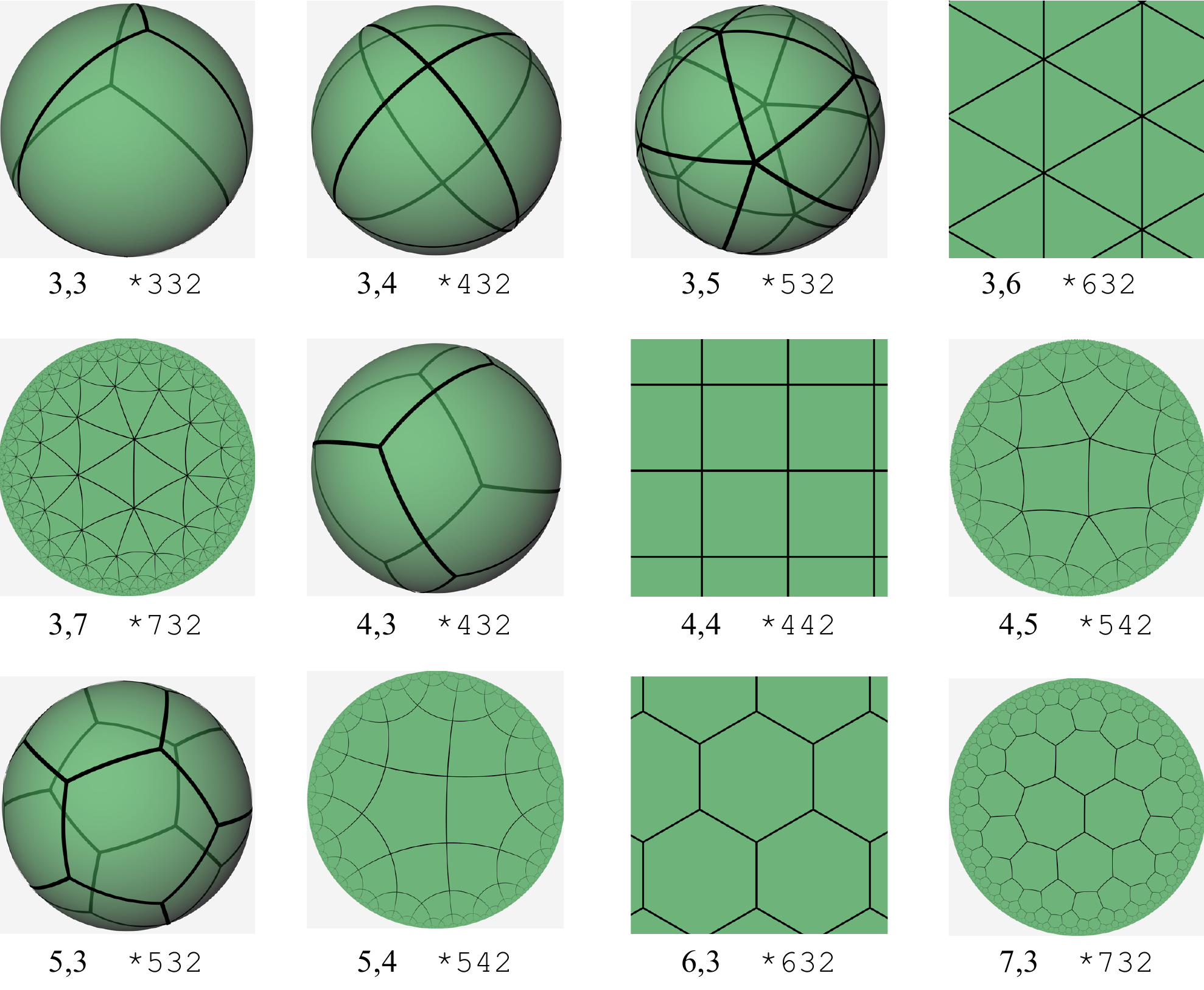}
    \caption{The 12 geometry-minimal types of periodic two-dimensional tilings with Dress complexity 1. Each labeled with $p,q$, that is, their tile and vertex 
    degress, and the orbifold name of their symmetry group.}

    \label{fig:size1}
\end{figure}

We can now formulate our first result: there exist exactly 12 different 
equivariant types of geometry minimal, periodic two-dimensional tilings
with Dress complexity 1, see Figure~\ref{fig:size1}. 
Note here we only consider tilings with tiles of degree 3 or more.
In addition, there are two types of geometry-minimal tilings
with digons, they have parameters ($p,q$) $2,3$ and $2,4$.

An enumeration of all possible equivariant types of geometry minimal, periodic two-dimensional tilings, up to a given maximal Dress complexity
$D$, can be obtained by systematically enumerating all non-isomorphic
Delaney-Dress symbols of size $\leq D$ that are geometry minimal.

Our second result is: there exist exactly 50 
equivariant types of geometry minimal, periodic two-dimensional tilings with Dress complexity 2.

More generally, we have solved this classification up to Dress complexity $D=24$ and obtain the following main result:
There exist exactly $2,395,220,319$ equivariant types of geometry minimal, periodic two-dimensional tilings with Dress complexity $\leq 24$.
Summary statistics are provided in Table~\ref{tab:stats}.
We have generated and saved all Delaney-Dress symbols for these tilings,
and make them available as described below.

\begin{table}[]
    \centering
    \begin{tabular}{r|rrrr}
   $\delta$ & $\#$ Spherical &  $\#$ Euclidean &  $\#$ Hyperbolic & Total \\
    \hline
 1 & 	           5 & 	           3 & 	           4 &   	12 \\
 2 & 	          13 & 	          15 & 	          22 &   	50 \\
 3 & 	          15 & 	           8 & 	          13 &   	36 \\
 4 & 	          30 & 	          37 & 	          71 &   	138 \\
 5 & 	          26 & 	          15 & 	          41 &   	82 \\
 6 & 	         119 & 	          86 & 	         221 &   	426 \\
 7 & 	         104 & 	          64 & 	         201 &   	369 \\
 8 & 	         252 & 	         217 & 	         796 &   	1,265 \\
 9 & 	         296 & 	         185 & 	         858 &   	1,339 \\
10 & 	         697 & 	         527 & 	        2,974 &   	4,198 \\
11 & 	         771 & 	         506 & 	        3,993 &   	5,270 \\
12 & 	        2,014 & 	        1,573 & 	       13,987 &   	17,574 \\
13 & 	        2,364 & 	        1,575 & 	       22,162 &   	26,101 \\
14 & 	        5,428 & 	        4,227 & 	       75,270 &   	84,925 \\
15 & 	        6,627 & 	        4,528 & 	      140,024 &   	151,179 \\
16 & 	       15,103 & 	       12,078 & 	      475,445 &   	502,626 \\
17 & 	       18,622 & 	       13,105 & 	      982,726 &   	1,014,453 \\
18 & 	       42,881 & 	       34,242 & 	     3,327,350 &   	3,404,473 \\
19 & 	       53,588 & 	       38,470 & 	     7,419,771 &   	7,511,829 \\
20 & 	      120,496 & 	       98,076 & 	    25,029,758 &   	25,248,330 \\
21 & 	      151,234 & 	      111,145 & 	    58,815,127 &   	59,077,506 \\
22 & 	      340,744 & 	      280,574 & 	   197,482,678 &   	198,103,996 \\
23 & 	      428,769 & 	      322,102 & 	   482,898,722 &   	483,649,593 \\
24 & 	      965,620 & 	      805,130 & 	  1,614,643,799 &   	1,616,414,549 \\
\hline
Total &   2,155,818    & 1,728,488 & 2,391,336,013 & 2,395,220,319\\
    \end{tabular}
    \caption{For Dress-complexity $\delta=1,\dots,24$, we list the number
    of different geometry-minimal periodic tilings of the sphere, euclidean plane and hyperbolic plane.}
    \label{tab:stats}
\end{table}

%%%%%%%%%%%%%%%%%%%%%%%%%%%%%%%%
\section{Visualization}

In Table~\ref{tab:stats} we count billions of Delaney-Dress symbols
that correspond to two-dimensional periodic tilings. 
To enable the exploration of these, we require an algorithm for calculating a drawing of the tiling associated with any given Delaney-Dress symbol $(\D,m)$.

Figure~\ref{fig:delaney} illustrates
that each node of a Delaney-Dress symbol corresponds to a different equivalence class of triangles in the barycentric subdivision of the corresponding periodic tiling, and that a fundamental domain for the
symmetry group can be obtained by selecting a suitable set of representatives
of the different classes of triangles.

To construct a tiling associated with a given two-dimensional Delaney-Dress symbol $(\D,m)$, we  proceed in three stages:

\begin{enumerate}
\item Compute coordinates for a barycentric triangulation
of the tiling for a fundamental domain of the symmetry group.
\item 
Compute a set of isometric transformations that generate the symmetry group.
\item Apply the generators to copies of
the triangulation of the fundamental domain so as cover
a desired region of the tiling.
\end{enumerate}

Stage (1) uses an algorithm and code developed by Klaus Westphal \citep{Westphal1991}.
The algorithm assigns a triangle to each node or chamber of the given
Delaney-Dress symbol $(\D,m)$. Triangles of adjacent nodes are then identified along the corresponding side, in a iterative manner.
This is done is such a way that the resulting triangulated region is a topological disc $B$ and all vertices of the triangulation that
correspond to a cone or corner point are located on the boundary.
Let $Z$ be a  $i,j$-component. 
We use $s(Z)$ to denote
the number of vertices of the triangulation that are associated with $Z$.
This will be one, if $Z$ lies in the interior of $B$.
If $s(z)>1$, then the vertices must all lie on the boundary of $Z$ and we say that $Z$ is {\em split}. 
For example, in Figure~\ref{fig:delaney}, the $1,2$-component containing chambers $1,2,3,5$ is split and is represented {\em twice} by vertices on the boundary of the fundamental domain, once involving the chambers labeled $1-3$ and the other time involving $5$.

Taking splitting into account, we assign an interior angle to $Z$ as
$\alpha(Z)=\frac{360^\circ}{s(Z)\times v_{ij}(D)}$, if $Z$ is a $i,j$-cycle, and
$\alpha(Z)=\frac{180^\circ}{s(Z)\times v_{ij}(D)}$, otherwise.

With this, we setup a polygon whose corners are given by the vertices assigned to the boundary of $B$, using the calculated interior angles. The polygon is heuristically fitted around an incircle and vertices with interior angle $180^\circ$ are then placed equally-spaced along the sides of the polygon.
Then all triangulation vertices that are assigned to the interior of $B$ are
iteratively assigned to the centroid of all adjacent vertices so as to obtain useful coordinates.

To address stage (2), a set of generators for the symmetry group is obtained
as follows. For a given chamber $D$ whose $i$-th edge lies on the boundary of $B$,
let $D'$ be its $i$-neighbor. Then a generator of the symmetry group can be obtained by calculating the
isometry that maps that the $0$-, $1$- and $2$-vertices of
the triangle representing $D$ onto the $0$-, $1$- and $2$-vertices of
the triangle representing $D'$, respectively. This is performed on all boundary chambers.

In stage (3), we repeatedly concatenate generators and keep
the transformed copy of the fundamental domain, if it will be visible.
The key practical challenge is to avoid placing more than one copy of the fundamental domain at the same location. To address this, we select a reference point within the interior of the fundamental domain and use either a quad-tree (in the case of euclidean tiles), or
a oct-tree (for spherical and hyperbolic tilings), to determine whether
 the current transformation applied to the reference point gives rise to a point that has already been seen.
 
 All hyperbolic calculations are performed using the Minkowski hyperboloid model. Visualization of the Poincare model and Klein model are implemented by observing the hyperboloid model using a perspective
 camera at locations $(0,0,-1)$ and $(0,0,0)$, respectively.
 
 To speed-up the visualization of translations, copies of the fundamental domain that disappear from view on one side of the tiling are reused and reappear on the other side of the tiling. 

%%%%%%%%%%%%%%%%%%%%%%%%%%%%%%%
\section{Enumeration and visualization software}

We have implemented the enumeration of Delaney-Dress symbols in a program called
{\em genDSyms} using the programming language Julia \citep{bezanson2017julia}.
For performance purposes, we use the principle of {\em orderly
generation} \citep{read1978every} to ensures that every symbol is produced exactly
once and no additional effort is required to identify and remove duplicates.
The process has two stages.

In the first stage we enumerate all possible Delaney-Dress {\em graphs} up to the required size.
Note that, for any given Delaney-Dress graph and choice of an initial node,
there exists a unique {\em ordered traversal},
that is, a breadth-first graph traversal, in which at each node
the incident edges are visited in the order of their labels \citep{delgado2003data}.
We use this to assign numbers to the nodes in the order they are encountered in and represent the traversal as a linear string of numbers by
listing the 0-, 1- and 2-neighbors of all the vertices in that same order.

As an example, consider the Delaney-Dress graph 
in Figure~\ref{fig:delaney}(b).  Beginning at the node on the left labelled
$1$, we see that it is its own 0- and 1-neighbor.  Its 2-neighbor is thus labelled $2$, its
1-neighbor in turn $3$, and so on.  Continuing in this fashion, we see that in fact the vertices are
already numbered in accordance with this traversal, which is then represented by the list
$1,1,2;\, 2,3,1;\, 4,2,5;\, 3,6,7;\, 7,5,3;\, 6,4,8;\, 5,8,4;\, 8,7,6$.
Of all the possible
traversals for this graph, this one turns out to be the lexicographically
smallest, because the leftmost node is the only one with both a 0- and a 1-loop, and thus the
only one that can give rise to a traversal representation starting with two ones.

To perform an orderly generation of Delaney-Dress graphs up to a given
size, we generate all possible ordered traversals and keep only those that are lexicographically
smallest for the graph they represent.
To speed up the enumeration process, we prune the enumeration tree 
whenever we identify a partial ordered traversal that cannot be completed to a 
lexicographically smallest one.

In the second stage of the enumeration,
for each Delaney-Dress graph we generate all possible valid definitions of the maps $m_{01}$ and $m_{12}$, in particular making use of the restrictions imposed by geometric minimality.

The file containing the complete galaxy of tilings is $322$~GB in size,
and it is thus impractical to make the file available on a webserver.
To provide easy access to much of the classification, we have
produced three SQLITE databases of Delaney-Dress symbols.
The first, {\tt tilings-1-18.tdb}, contains all tilings
with Dress complexity 1--18. The other two, {\tt spherical-1-24.tdb}
and {\tt euclidean-1-24.tdb},
contain all spherical and euclidean tilings, respectively of Dress complexity 1--24.
Each database contains a table called ``tilings'' that has the schema
shown in Table~\ref{tab:schema}.

\begin{table}[]
    \centering
    \begin{tabular}{lp{8cm}}
    Column name and type & Explanation\\
\hline
        id INTEGER PRIMARY KEY & number in file\\
 symbol TEXT & Delaney-Dress symbol $(\D,m)$\\
 complexity INTEGER & Dress complexity $\delta(\D,m)$\\
 geometry TEXT &  Name of two-dimensional geometry\\
 curvature TEXT & Curvature ${\mathcal K}(\D,m)$\\
 euler REAL & Euler characteristic\\ 
 orbifold TEXT &  Orbifold name of symmetry group\\
 symmetry\_class TEXT & Symmetry class of graph, as defined in \citep{Hyde:2014aa}\\
 signature TEXT &  Signature is expression such as $(3 4 6 5)$ that indicates the tiling consists of tiles of degree 4 with vertices of degree
 3, 4, 6 and 5.\\
  tile\_deg TEXT & List of tile degrees in ascending order\\
 vertex\_deg TEXT & List of vertex degrees in ascending order\\
 tiles INTEGER & Number of equivalence classes of tiles\\ 
 edges INTEGER & Number of equivalence classes of edges\\
 vertices INTEGER & Number of equivalence classes of vertices\\
 normal BOOLEAN & Is tiling pseudo-convex?\\
 maximal BOOLEAN &  Is the symmetry group of the tiling maximal? \\
 colorable BOOLEAN &  Tiling is colorable, if no two tiles of the same symmetry equivalence class share an edge. \\
  orientable BOOLEAN &  Does symmetry group only contain orientation-preserving symmetries?\\
 fixed\_point\_free BOOLEAN &  Is symmetry group fixed-point free?\\
  self\_dual BOOLEAN & Is tiling self-dual?\\
    \end{tabular}
    \caption{Schema for table ``tilings'' used for storing Delaney-Dress symbols and some associated properties.}
    \label{tab:schema}
\end{table}

We have implemented a new program called Tegula that can be used to explore
our galaxy of tilings. Tegula takes as input an SQLITE database of
Delaney-Dress symbols and provides drawings of the corresponding tilings
in a ``collection'' tab, on a page-by-page basis.
The program provides an interactive dialog for searching
for tilings of specific interest. The user can page through all
tilings that fulfill a given query. 

\begin{figure}
    \centering
    \includegraphics[width=\textwidth]{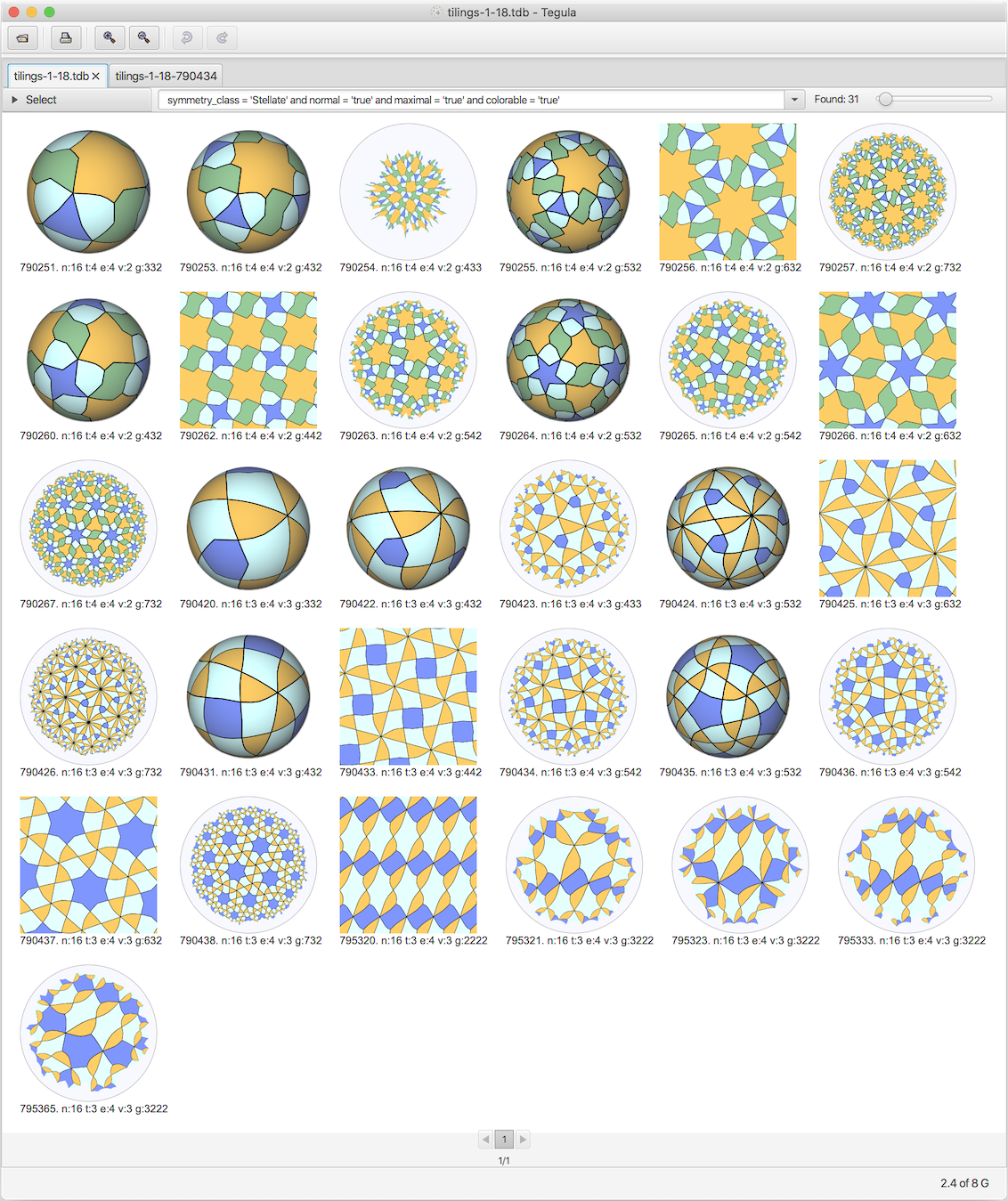}
    \caption{Applying the query  {\tt symmetry\_class = 'Stellate' and normal = 'true' and maximal = 'true' and colorable = 'true'} to the database {\tt tilings-1-18.tdb} of
    all Delaney-Dress symbols of Dress complexity $\leq 18$
    returns 31 tilings.
    }
    \label{fig:db-example}
\end{figure}

For example,  the query  ``symmetry\_class = 'Stellate' and normal = 'true' and maximal = 'true' and colorable = 'true' '' returns 31 tilings
of Dress complexity $\leq 18$, 
displayed Figure~\ref{fig:db-example}.

\begin{figure}
    \centering
    \includegraphics[width=\textwidth]{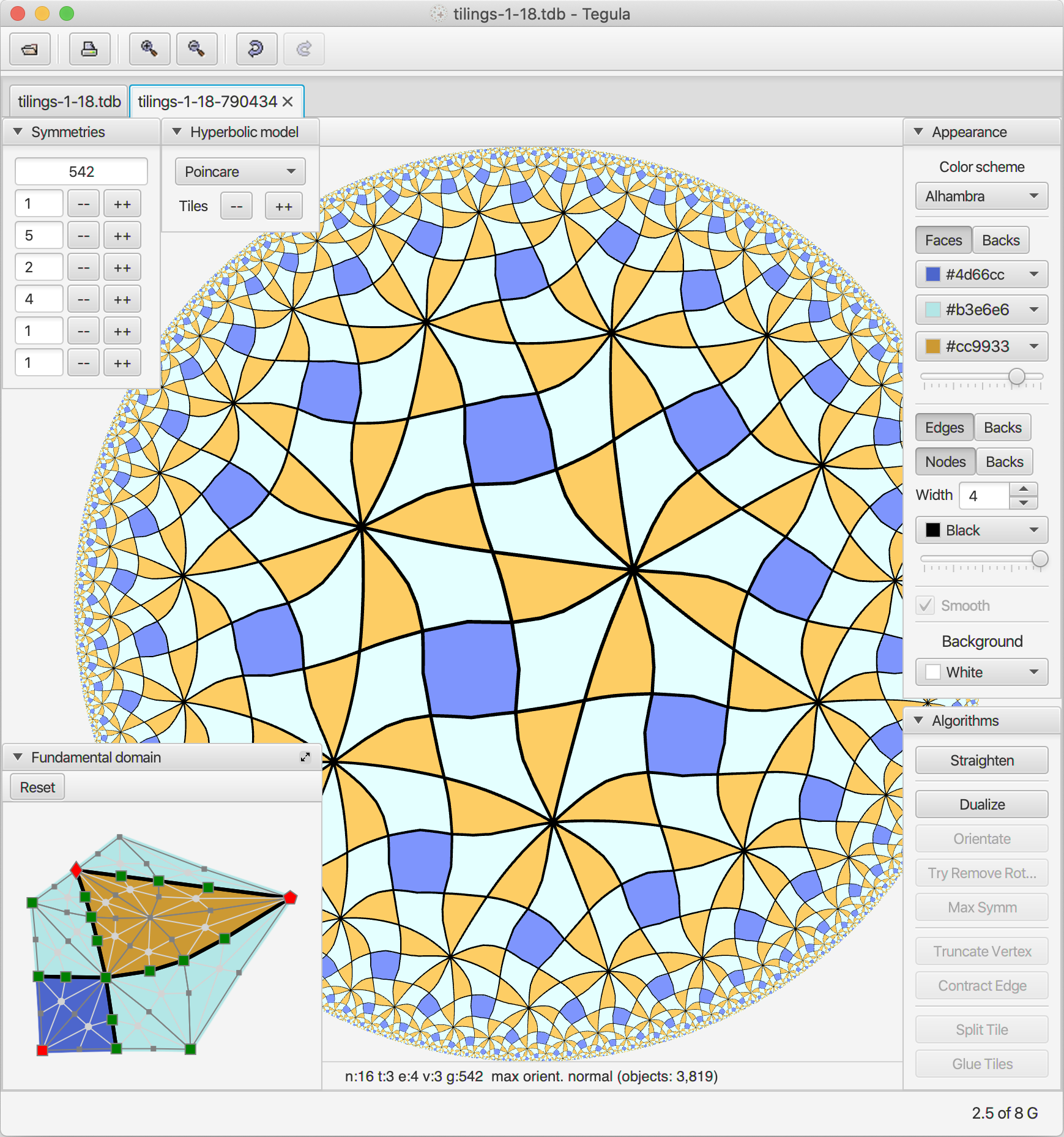}
    \caption{
    Any tiling can be edited using five different panels.
    The {\em symmetries} panel allows editing of the rotational orders
    of the symmetry group; the {\em hyperbolic model} panel allows
    selection between the Poincare, Klein and hyperboloid
    model; the {\em appearance} panel allows the modification 
    of different aspects of the tiling; the {\em algorithms}
    panel provides some transformations of the drawing and
    the Delaney symbol; and the {\em fundamental domain} panel
    allows interactive reshaping of the edges of the tiling.
    }
    \label{fig:editor}
\end{figure}

Individual tilings can be edited in a number of different ways.
Double-clicking on a tiling in a collection tab will open the tiling
in a new ``editor'' tab and there five panels of tools are available to
modify the displayed tiling, as illustrated in Figure~\ref{fig:editor}.

\section*{Availability}

The enumeration program genDSyms is written in Julia. The source is provided here:
\url{ https://github.com/odf/julia-dsymbols}.
The visualization and exploration program Tegula is written in Java and uses the OpenJFX library.
The source is provided here:
\url{https://github.com/husonlab/tegula}.
Installers for Windows and MacOS are available here:
\url{https://software-ab.informatik.uni-tuebingen.de/download/tegula}.

\section*{Acknowledgements}

We thank Klaus Westphal for providing us with his 
original source code for computing the fundamental domain of a tiling.
We thank Julius Vetter and Cornelius Wiehl for programming contributions.

\bibliographystyle{abbrvnat}
\bibliography{references}
\end{document}